\documentclass[11pt]{amsart}
\usepackage{geometry,amsrefs,hyperref}                
\usepackage{graphicx}
\usepackage{amssymb}
\usepackage{epstopdf}
\newtheorem{prop}{Proposition}
\newtheorem{theo}{Theorem}
\newcommand{\I}{\operatorname{i}}
\newcommand{\C}{\mathbb C}

\newcommand{\e}{\operatorname{e}}
\renewcommand{\Re}{\operatorname{Re}}
\renewcommand{\Im}{\operatorname{Im}}

\title{Explicit description of spherical rigid hypersurfaces in $\mathbb C^2$}
\author{Vladimir Ezhov \and Gerd Schmalz}

\begin{document}
\maketitle
\section{Introduction}

The class of real hypersurfaces in complex space that are invariant with respect to an infinitesimal translation transversal to the complex tangent space is known as rigid hypersurfaces (see \cite{BRT}). In this article we consider rigid hypersurfaces in $\mathbb C^2$ with coordinates $z,w=u+\I v$. In this case a rigid hypersurface can be locally described by an equation of the form $v=h(z)$. 

In 1991 Stanton \cite{Sta91} developed a normal form for rigid hypersurfaces. Other normal forms that reflect the presence of symmetries have been constructed by Kol\'a\v{r} \cite{Kol} and Ezhov, Kol\'a\v{r} and Schmalz \cite{EKS} more recently. 

It is a natural question how to recognise hypersurfaces that are holomorphically equivalent to the model hypersurface, i.e. the sphere $v=|z|^2$, via its normal form. As an application of the normal form Stanton derived a list of examples of normal forms equivalent to the sphere:

\begin{align}\label{sinh}
\frac{1}{2r}\sinh 2r v =& |z|^2\\\label{Stanton}
\frac{1}{2r}\sin 2rv \left(1-\frac{2|b|^2\theta}{|c|^2} \right)=&\frac{|z|^2 \e^{-2\theta v}}{1+4|b|^2|z|^2+ 2\I(b\bar{z}-\bar{b} z)} + \frac{|b|^2}{|c|^2}(\e^{-2\theta v}-\cos 2rv)+ \\\nonumber &+ \frac{\bar{b} z}{\bar{c}(1-2\I\bar{b}z)}(\e^{-2\theta v}-\e^{2\I rv})+ \frac{b\bar{z}}{c(1+2\I b\bar{z})}(\e^{-2\theta v}-\e^{-2\I rv}),
\end{align} 
where $b$ and $c=r+\I \theta$ are complex parameters. The remaining examples from Stanton's list can be obtained by letting $r$ and $\theta$ converge to $0$. In particular, for $b=0$ and $\theta=0$ one finds the pendant to (\ref{sinh})
$$\frac{1}{2r}\sin 2r v = |z|^2.$$

Stanton has raised the question whether this list was complete. Her method was based on the analysis of the holomorphic mappings that transform the infinitesimal translation  $\partial/\partial w$ into a vector field of a 7-parameter family of infinitesimal sphere automorphisms. She acknowledged that the listed examples correspond to a 4-parameter subfamily. 

The main aim of this paper is to provide a complete classification of spherical rigid hypersurfaces. Combined with Stanton's results this also answers Rothschild's question which rigid hypersurfaces are equivalent to a hypersurface that is given by an equation
$$v= p(z,\bar{z})$$ 
where $p$ is a homogeneous polynomial (see \cite{Sta91}). Stanton has solved this problem for polynomials of degree bigger than two. Our result completes the case of degree $2$.

The following argument gives the upper bound 4 for the number of parameters involved: In the Levi-nondgenerate case, which applies for spherical rigid hypersurfaces, Stanton's normal form is similar to Chern-Moser's normal form, except for the terms of the defining function $\phi$ of bidegrees (2,2), (2,3), (3,2) and (3,3) with respect to $z,\bar{z}$ being allowed to be any constants. These four real constants completely control the rigid hypersurface. Thus the normal form equation becomes
$$v= |z|^2 + c_{22}|z|^4 + c_{23}z^2\bar{z}^3 + c_{32}z^3\bar{z}^2 + c_{33}|z|^6 + \sum_{\stackrel{\min(j,k)\ge 2}{\max(j,k)\ge4}} F_{jk} z^j\bar{z}^k.$$

The coefficients $c_{22}$, $c_{23}=\bar{c}_{32}$ and $c_{33}$ depend in an algebraic way on Stanton's parameters $b$ and $c$. It can be shown that not all coefficients can be attained in this way, which indicates that Stanton's list is incomplete. In the case when $c_{23}=0$ Stanton's example reduces to
$$\frac{1}{2r}\sin 2rv=  \e^{-2\theta v}|z|^2,$$
which covers only the coefficients $c_{33}\ge \frac{3}{2}c_{22}^2$. The coefficients $c_{33}< \frac{3}{2}c_{22}^2$ are covered by the family
$$\frac{1}{2r}\sinh 2rv=  \e^{-2\theta v}|z|^2,$$
which for $\theta=0$ coincides with (\ref{sinh}) but was absent in Stanton's list for $\theta\neq 0$. 

We show that Stanton's mappings of a rigid sphere to the Heisenberg sphere can be modified by combining them with suitable sphere automorphisms so that all parameters can be covered. 

This yields not only the desired classification but also provides the complete solution to a non-linear PDE that expresses the zero-curvature equation of a rigid hypersurface $M$ in normal form 
\[v= h(z).\]
It is well known that local sphericity is equivalent to vanishing of the Cartan curvature, which was computed in \cite{EMS}. For rigid hypersurfaces the zero-curvature equation simplifies to
\begin{equation} \label{pde}
f_{z\bar{z}\bar{z}\bar{z}} -3 f_{z\bar{z} \bar{z}}f_{\bar{z}}  +2  f_{z\bar{z}}(f_{\bar{z}})^2 - f_{z\bar{z}} f_{\bar{z}\bar{z}} =0
\end{equation}
where $f(z)=\log\Delta h(z)$. \bigskip

Acknowledgements. The authors are grateful to Martin Kol\'a\v{r} for numerous useful discussions from which some of the ideas used in this article arose. The research was supported by the Max-Planck-Institut f\"ur Mathematik and the ARC Discovery grant DP130103485.\bigskip

\section{Modified normalisation}
We apply a modified Chern-Moser normalisation procedure to the Heisenberg sphere $\Im w_2=|z_2|^2$ to obtain a hypersurface of the form
\[v= |z|^2 + c_{22}|z|^4+ c_{23}z^2\bar{z}^3 + c_{32}z^3\bar{z}^2+c_{33}|z|^6+ \cdots\] 
For given $c_{22}, c_{23}, c_{33}$ such normalisation mapping is uniquely determined up to automorphisms of the Heisenberg sphere. All rigid spheres can be found by this procedure, though the resulting hypersurface does not have to be rigid in the higher order terms a priori. We construct the inverse mapping as a composition of two:
\begin{align*}
 z_2&= p(w_1) + z_1 + 2\I \sum_{j=2}^\infty T_j(w_1) z_1^j=p(w_1) + \frac{z_1}{1-2\I \bar{p}'(w_1)z_1}\\
 w_2&= q(w_1) + 2\I  \sum_{j=1}^\infty g_j(w_1) z_1^j= q(w_1) + \frac{2\I \bar{p}(w_1)z_1}{1-2\I \bar{p}'(w_1)z_1}
\end{align*}
Here we assume 
\begin{equation}\label{chainpar}
q'= 1+2\I p'\bar{p}.
\end{equation}
Condition \eqref{chainpar} is imposed for convenience in computation and can be fixed at the next step.

The second mapping has the form
\begin{align*}
 z&= \e^{\I \alpha (w_1)} \sqrt{h'(w_1)} z_1\\
 w&= h(w_1).
\end{align*}

Then $p,\alpha,h$ satisfy the equations
\begin{align}\label{norm1}
6 |p'|^2 +2  \alpha '-c_{22} h'&=0\\\label{norm2}
-c_{23} e^{-i \alpha} (h')^{3/2}- 2 p'' +4\I |p'|^2 p'&=0\\\label{norm3}
\frac{1}{3} h'''-\frac{(h'')^2}{2 h'} +\frac{3c_{22}^2-2c_{33}}{2} (h')^3 
+2 |p'|^4h' +\frac{2}{3}(\I p'' \bar{p}'-\I p'\bar{p}'') h'&=0
\end{align}
with initial conditions $p(0)=0$, $h(0)=0$, $h'(0)>0$. Using combinations with automorphisms of the Heisenberg sphere we may assume that $\alpha(0)=0$, $h'(0)=1$, $p'(0)=0$, $h''(0)=0$. It follows that, for given $c_{22}, c_{23}, c_{33}$, system\eqref{norm1}-\eqref{norm3} has unique solutions $\alpha, p, h$.

Though it is easy to successively compute the coefficients of the power series of $\alpha,p,h$ from recursive formulae we have not succeeded in solving this system of ODE, except for $c_{23}=0$. In this case $p\equiv 0$ and 
\begin{align*}
h(u)&= \frac{2}{\sqrt{9 c_{22}^2- 6c_{33}}} \arctan \frac{\sqrt{9 c_{22}^2- 6c_{33}}u}{2}\\
\alpha(u)&= \frac{c_{22}}{2} h(u).
\end{align*}

\section{Stanton's mapping}
Rigid surfaces can be characterised by the presence of a translation symmetry transversal to the complex tangent spaces. In suitable local coordinates such symmetry is generated as the flow of the vector field $\frac{\partial}{\partial u}$. Stanton constructed holomorphic mappings in the ambient space that pull a suitable infinitesimal automorphism $X$ of the Heisenberg sphere back to  $\frac{\partial}{\partial u}$. The relevant infinitesimal automorphisms of the Heisenberg sphere are well-known and form a 7-parametric family consisting of
\[X= 2\Re \left((b+ cz+aw+2\I \bar{a}z^2+ \rho zw)\frac{\partial}{\partial z} + (1+2\I\bar{b}z+ 2rw+2\I \bar{a}zw+\rho w^2)\frac{\partial}{\partial w} \right),\]
where $b,a\in\C$ and $c=r+\I\theta\in \C^*$ and $\rho\in\mathbb R$.

The resulting system is   
\begin{align}\label{system}
\frac{\partial Z}{\partial w} &= b+ cZ+aW+2\I \bar{a}Z^2+ \rho ZW\\\nonumber
\frac{\partial W}{\partial w} &= 1+ 2\I \bar{b} Z+ 2r W+2\I \bar{a}ZW+\rho W^2
\end{align}
with initial conditions $Z(z,0)\equiv \frac{z}{1-2\I \bar{b}z}$, $W(z,0)\equiv 0$. Stanton solved this system in the particular case $a=\rho=0$, that is when \eqref{system} is linear.
The solutions are
\begin{align}\label{stmap}
Z&= \frac{b}{c}(\e^{cw}-1) +\frac{\e^{cw}z}{1-2\I \bar{b}z}\\\nonumber
W&= \left(1-\frac{2\I|b|^2}{c}\right)\frac{\e^{2rw}-1}{2r} + \frac{2\I\bar{b}}{\bar{c}}\left(\frac{z}{1-2\I\bar{b}z}+\frac{b}{c}\right)(\e^{2rw}-\e^{cw}).
\end{align}
Substituting this mapping into the sphere equation $v=|z|^2$ yields (\ref{Stanton}). We show that not all values of the parameters $c_{22},c_{23},c_{33}$ can be realised by these mappings. Direct computation shows that 
the explicit 6$^{\text{th}}$ order expansion of (\ref{Stanton}) is 
\begin{multline} \label{series}
 v= |z|^2 + (6|b|^2-2\theta) |z|^4 +(2 \bar{c} +4 \I |b|^2) b z^2\bar{z}^3 +(2 c -4 \I |b|^2) \bar{b} z^3\bar{z}^2 + \\ + \left(\frac{2}{3}r^2 + 6 \theta ^2 +56 |b|^4 -\frac{112}{3} \theta  |b|^2\right) |z|^6+ \cdots\end{multline}
Hence the parameters take the form
\begin{align}\label{para22}
c_{22}&= 6 |b|^2 -2\theta\\\label{para23}
c_{23} &= 2 (r-\I \theta )b +4 \I b |b|^2\\\label{para33}
c_{33} &= \frac{2}{3}r^2 + 6 \theta ^2 +56 |b|^4 -\frac{112}{3} \theta  |b|^2.
\end{align}
Stanton's family realises all those $c_{22},c_{23},c_{33}$ for which the system of algebraic equations \eqref{para22}-\eqref{para33} has a solution $b\in\mathbb C$, $r,\theta\in\mathbb R$. A method to solve this system is to first express $\theta$ using equation \eqref{para22} then express $r^2$ through $c_{33}$ and $|b|^2$ using equation \eqref{para33} and plug it into the absolute square of equation \eqref{para23}. Solve the resulting cubic equation on $|b|^2$ and finally find $r,\theta$. However the resulting values for $|b|^2$ or $r^2$ could turn out to be negative.\medskip

Consider two examples: 1. Let $c_{23}=0$. Then $b=0$ and \eqref{para33} becomes 
$$r^2=\frac{3}{2}c_{33}- \frac9{4} c_{22}^2.$$
By allowing $r^2$ to become negative the solution
$$\frac{1}{2r}\sin 2rv=e^{-2\theta v}|z|^2$$
turns into \eqref{sinh}.

2. For $c_{22}=0$ and $c_{23}=2$ only $c_{33}\ge -2$ is feasible. Indeed, from the equation (\ref{para22}) we get $\theta=3 |b|^2$. Then (\ref{para23}) is equivalent to
$$\frac1b= r - \I |b|^2.$$
This implies
$$r^2 + |b|^4=\frac{1}{|b|^2}$$
and hence
$$r^2= \frac{1}{|b|^2}- |b|^4.$$
This is only possible if $|b|\le 1$.
Now (\ref{para33}) becomes  
$$c_{33}=\frac23\frac{1}{|b|^2} -\frac83 |b|^4\ge -2.$$  

Finally, we point out that Stanton's mapping can be interpreted as modified normalisation with
\begin{align*}
\alpha(u)&= \frac{-\theta}{2r} \log(1+2ru)\\
h(u)&=\frac{1}{2r} \log(1+2ru)\\
p(u)&= \frac{b}{c} \left(\e^{\frac{c}{2r}\log (1+2ru)} -1  \right).
\end{align*}

Notice that the initial conditions are $p'(0)=b$ and $h''(0)=-2r$.\medskip

On the other hand, the proposition below shows that the solutions of \eqref{system} with different vector fields $X$ yield the missing rigid spheres.

\begin{prop}
For any set of parameters $c_{22},c_{23},c_{33}$ there exists a rigid sphere. It can be realised by
the solutions $Z(z,w), W(z,w)$ of \eqref{system} with parameters $c=\I\theta, a,\rho$ and $\Re c=b=0$, where 
\begin{align*}
\theta&=-\frac{c_{22}}{2} \\
a&= -\frac{c_{23}}{2}\\
\rho&=  -\frac32 c_{33} + \frac94 c_{22}^2
\end{align*}
and initial conditions $Z(z,0)=z$, $W(z,0)=0$.
\end{prop}

Proof. For any choice of parameters $c_{22},c_{23},c_{33}$ we find the corresponding parameters $\theta, a, \rho$ and hence the infinitesimal sphere automorphism $X$. According to the Picard-Lindel\"of theorem the system \eqref{system} has a unique solution with initial conditions $Z(z,0)=z$, $W(z,0)=0$.

Direct computations with jets show that
\begin{align*}
Z(z,w)&=z+\I\theta zw +\frac{a w^2}{2}+2\I \bar{a} z^2 w +\frac{\rho -\theta^2}{2} zw^2  + \frac{\I\theta a}{6} w^3+\cdots\\
W(z,w)&= w+ \I \bar{a} z w^2 +\frac{\rho  w^3}{3}+\cdots
\end{align*}
where the dots indicate terms of order higher than 3. Substituting these truncated mappings into the sphere equation 
$$\frac{W-\bar{W}}{2\I}-Z\bar{Z}=0$$
yields  
$$v= |z|^2 - 2\theta |z|^4 -2a z^2\bar{z}^3-2 \bar{a} z^3\bar{z}^2 + \left( 6\theta^2- \frac{2\rho}{3}\right) |z|^3+\cdots$$
where the dots indicate terms of bidegree (2,4) and (4,2) and higher order. 
\hfill$\Box$\medskip

It is an immediate consequence that the modified Chern-Moser normalisation of the Heisenberg sphere from Section 2 yields indeed rigid hypersurfaces.

\section{Twisting Stanton's mapping}
One way to produce solutions of \eqref{system} with non-trivial $a,\rho$ is to apply the adjoint action of $SU(2,1)$ on the coefficient matrix of the linear system in a suitable way. Geometrically this amounts to compose Stanton's mapping with a sphere automorphism. The composition of Stanton's mapping $z_1(z,w), w_1(z,w)$ with parameters $c=r+\I \theta, b$ and the sphere automorphism
\begin{align*}
z_2&=\frac{z_1-b w_1}{1+2\I \bar{b} z_1 +(r-\I |b|^2)w_1}\\
w_2&=\frac{w_1}{1+2\I \bar{b} z_1 +(r-\I |b|^2)w_1}
\end{align*}   
is
\[z_2=\frac{P_1}{Q},\qquad
w_2=\frac{P_2}{Q}\]
with
\begin{align}\label{mapping}
P_1=& (a+4(\theta-\phi)\phi z )\frac{\cosh rw - \e^{\I\theta w} + \I\theta \frac{\sinh rw}{r}}{r^2+\theta^2} +\left(-2\I\phi\frac{\sinh rw}{r}+ \e^{\I\theta w}\right)z \\\nonumber
P_2=&2\I(\phi-\bar{a} z)\frac{\cosh rw - \e^{\I\theta w}+ \I\theta \frac{\sinh rw}{r}}{r^2+\theta^2}+ \frac{\sinh rw}{r}  \\\nonumber
Q=& 2 (\phi-\theta)(\phi-\bar{a}z)\frac{\cosh rw - \e^{\I\theta w}+ \I\theta \frac{\sinh rw}{r}}{(r^2+\theta^2)} +\I(\phi-2\bar{a} z) \frac{\sinh rw}{r}+\cosh rw 
\end{align}
where $a = -b (r-\I \theta +2 \I |b|^2)$ and $\phi=|b|^2$. 

Notice that the Taylor series of 
\[\cosh rw -\e^{\I\theta w} + \I\theta \frac{\sinh rw}{r}\] 
has coefficients $a_0=a_1=0$
$$a_{2n+1}= \frac{\I\theta(r^{2n}-(-1)^n \theta^{2n})}{(2n+1)!},\qquad a_{2n}=\frac{r^{2n}-(-1)^n \theta^{2n}}{(2n)!}$$
for $n\ge 1$ and therefore is divisible by $r^2+\theta^2$. Clearly $\sinh rw$ is divisible by $r$. Therefore, $P_1,P_2,Q$ are entire functions with respect to $z,w,r,\theta,a,\phi$. Moreover they are even functions with respect to $r$.

The vector field $\frac{\partial}{\partial w}$ is pulled back to
\[( \I \tau z_2 + aw_2 +2\I \bar{a}z_2^2 + \rho z_2w_2)\frac{\partial}{\partial z_2} + (1+ 2\I\bar{a}z_2w_2 + \rho w_2^2 )\frac{\partial}{\partial w_2},\]
where
\begin{align}\label{alg}
\tau&=\theta -3\phi=-\frac{c_{22}}{2} \\\nonumber
a&= -b (r-\I\theta + 2\I \phi)= -\frac{c_{23}}{2}\\\nonumber
\rho&= -3 \phi^2 - r^2 + 2 \phi \theta= -\frac32 c_{33} + \frac94 c_{22}^2.
\end{align}

We make formula \eqref{mapping} universal by allowing imaginary $r$ and negative $\phi$.

Solving \eqref{alg} for $\theta,b,r$ we find
\begin{align}\label{paratheta}
\theta=& \tau + 3\phi\\\nonumber
r^2=& -\rho + (2\tau+3\phi)\phi
\end{align}
Now, $\phi$ can be determined from the equation 
\begin{equation}\label{algset}
|a|^2= 4\phi^3+ 4\tau\phi^2+ (\tau^2-\rho)\phi.
\end{equation}
Notice that the cubic equation has real coefficients and therefore has at least one real solution. Let $\phi$ be any real solution.

For $\tau=0$ the solution $\phi$ can be given by the formula
\[\phi=\frac12\left(\left(|a|^2 +\sqrt{|a|^4-\frac{\rho^3}{27}}\right)^{\frac13} + \left(|a|^2 -\sqrt{|a|^4-\frac{\rho^3}{27}}\right)^{\frac13} \right).\]
Here we take the principal branch of the cubic root on the right half-plane and the real cubic root on the real axis.
Notice that $\phi$ is a continuous real-valued function, though not necessarily non-negative, and that $\phi=0$ for $a=0$.

\begin{theo}\label{theo} Let 
\begin{align*}
Z(z,w,a,\rho,\tau,\phi)=\frac{P_1}{Q}\\
W(z,w,a,\rho,\tau,\phi)=\frac{P_2}{Q}
\end{align*}
where $P_1,P_2,Q$ are as in \eqref{mapping} and $\theta$ and $r^2$ are expressed as functions of $\tau, \rho, \phi$ by \eqref{paratheta}.

Then $Z,W$ satisfy the system \eqref{system} on the real algebraic set given by \eqref{algset}. 

It follows that all rigid spheres can be found as inverse images of $\Im W=|Z|^2$ under the mappings $Z,W$ with suitable parameters $\tau,a,\rho,\theta$ and therefore have the form
\begin{equation*}
(1+2\phi |z|^2) \frac{\sin 2rv}{2r} - \e^{-2\theta v}|z|^2-(\phi+ \bar{a}z+a \bar{z}+4\phi(\phi-\theta)|z|^2) \frac{\e^{-2\theta v}-\cos 2rv + \frac{\theta\sin 2rv}{r}}{r^2+\theta^2} =0.
\end{equation*}
\end{theo}

\noindent {\sc Proof.} Direct computation shows that the expressions
\begin{align*}
\frac{\partial Z}{\partial w} &- \I \tau Z-aW-2\I \bar{a}Z^2- \rho ZW\\\nonumber
\frac{\partial W}{\partial w} &- 1- 2\I \bar{a}ZW-\rho W^2
\end{align*}
factorise  with a factor $4\phi^3+ 4\tau\phi^2+ (\tau^2-\rho)\phi - |a|^2= 4 \phi ^3-4 \theta  \phi ^2 +(\theta ^2 +r^2) \phi -|a|^2$. 

In fact
\begin{multline*}
Q \frac{\partial P_1}{\partial w}-P_1 \frac{\partial Q}{\partial w}
-a P_2 Q-2 \I \bar{a} P_1^2+ \left( r^2+3 \phi ^2 -2 \theta  \phi\right) P_1 P_2-\I(\theta -3 \phi ) P_1 Q =\\
\frac{2 \left(4 \phi ^3-4 \theta  \phi ^2 +(\theta ^2 +r^2) \phi -|a|^2\right) \left(e^{\I \theta  w}- \cosh (r w)-\I \theta  \frac{\sinh (r w)}{r}\right)}{r^2+ \theta^2}\times \\
\Big[-\I  \left(2 \bar{a}  (\theta -\phi )z^2+(\theta ^2-6 \theta  \phi +r^2+6 \phi ^2)z-a\right)(e^{\I \theta  w}+\cosh (r w))\\+  \left(2 \bar{a}\left(\theta  \phi +r^2\right) z^2+ 3  \phi (\theta^2 -2 \phi\theta-r^2 )z \phi +\theta a\right)\frac{\sinh (r w)}{r}\Big]
\end{multline*}
and
\begin{multline*}
Q \frac{\partial P_2}{\partial w}-P_2 \frac{\partial Q}{\partial w}
- Q^2 -2 \I \bar{a} P_1P_2+ \left(r^2+3 \phi ^2 -2 \theta  \phi\right) P_2^2 =\\
\frac{2 \left(4 \phi ^3-4 \theta  \phi ^2 +(\theta ^2 +r^2) \phi -|a|^2\right) \left(e^{\I \theta  w}- \cosh (r w)-\I \theta  \frac{\sinh (r w)}{r}\right)}{r^2+ \theta^2}\times \\
\Big[2  (\bar{a} z-\phi ) (e^{i \theta  w} -\cosh (r w))-i \left( 2 \bar{a}\theta  z+\theta^2 -2 \phi \theta  +r^2\right)\frac{\sinh (r w) }{r}\Big].
\end{multline*}

Therefore $Z,W$ satisfy the system \eqref{system} on the real algebraic set \eqref{algset}. 

Now the rigid sphere formula can be obtained either from Stanton's formula \eqref{Stanton} by replacing $r,\theta,b$ by their expressions in $\tau,a,\rho,\theta$ or by inserting $Z,W$ into the standard Heisenberg sphere equation. 
 $\Box$\bigskip

Notice that the Taylor series of $\e^{-2\theta v}-\cos 2rv + \frac{\theta\sin 2rv}{r}$ has coefficients $a_0=a_1=0$
$$a_{2n+1}= -\frac{2^{2n+1}\theta(r^{2n}-(-1)^n \theta^{2n})}{(2n+1)!},\qquad a_{2n}=\frac{2^{2n}(r^{2n}-(-1)^n \theta^{2n})}{(2n)!}$$
for $n\ge 1$ and therefore is divisible by $r^2+\theta^2$. Clearly $\sin 2rv$ is divisible by $r$. It follows that the rigid sphere formula is an entire function with respect to all variables and parameters.\medskip

We find an example of a rigid sphere that is not in Stanton's family by setting
\[ \tau=0,\quad a=\sqrt{2},\quad \rho= 6.\]
Then 
\[\theta=-3, \quad \phi= -1, \quad r^2=-3\]
and 
\begin{equation*}
(1-2|z|^2) \frac{\sinh 2\sqrt{3}v}{2\sqrt{3}} - \e^{6 v}|z|^2+(1- \sqrt{2} (z+ \bar{z})+8|z|^2) \frac{\e^{6 v}-\cosh 2\sqrt{3} v -\sqrt{3} \sinh 2\sqrt{3}v}{6} =0.
\end{equation*}

\section{The zero-curvature equation}

Local equivalence of a real hypersurface $M$ in $\mathbb C^2$ to a sphere can be characterised by vanishing of its Cartan curvature. In \cite{EMS} an explicit expression of the Cartan curvature has been computed. In the case of rigid hypersurfaces this expression considerably simplifies. If $M$ is given by the equation
\[v=h(z) \]
vanishing of the Cartan curvature is equivalent to the non-linear PDE
\[ f_{z\bar{z}\bar{z}\bar{z}} -3 f_{z\bar{z} \bar{z}}f_{\bar{z}}  +2  f_{z\bar{z}}(f_{\bar{z}})^2 - f_{z\bar{z}} f_{\bar{z}\bar{z}} =0\]
with $f=\log \Delta h$. It reduces to
\begin{equation}\label{reducedpde} 
\tilde{f}_{z\bar{z}\bar{z}} -3 \tilde{f}_{z\bar{z}}\tilde{f}  +2  \tilde{f}_{z}\tilde{f}^2 - \tilde{f}_{z} \tilde{f}_{\bar{z}} =0
\end{equation}
with $\tilde{f}=\frac{\partial}{\partial \bar{z}}\log \Delta h$.

Notice that the normal form conditions on $h$ are encoded in the choice of the solution of the auxiliary equation
\[\frac{\partial}{\partial \bar{z}}\log \Delta h=\tilde{f}.\]
It follows that the implicit equation of the rigid spheres given in Theorem \ref{theo} provides the solutions to \eqref{reducedpde}.

Attempts to use the zero-curvature equation for solving the rigid sphere problem were unsuccessful. We only succeeded in solving special cases of the zero-curvature equation with the additional assumption that $M$ is  circular, i.e. $v=f(|z|^2)$,  or that $M$ is a tube, i.e. $v=h(x)$.  

In the circular case 
\[v= h(|z|^2)= |z|^2 + a_2 |z|^4 + a_3 |z|^6 +\dots\]
and the equation on $g=\log \Delta h$ with $t$ replacing $|z|^2$ is
\[
tg^{IV}+3 g'''-  g' \left(3 t g'''+7 g''\right)-t \left(g''\right){}^2+2 t \left(g'\right){}^2 g''+2 \left(g'\right){}^3=0.
\]

Formal power series solutions 
\[g(t)= c_1 t + c_2 t^2 +\dots \]
are determined by $c_1,c_2$. They correspond to the surfaces

\[\frac{\sin \alpha v}{\alpha}= \e^{-2\beta v} |z|^2\]

with $c_1=-4\beta$, $c_2=-\frac{3}{2}\alpha^2- \frac{43}{2}\beta^2$ and

\[\frac{\sinh \alpha v}{\alpha}=  \e^{-2\beta v} |z|^2\]
with $c_1=-4\beta$, $c_2=-\frac{3}{2}\alpha^2 + \frac{11}{2}\beta^2$.\bigskip

The spherical tubes are well-known (see, for instance, \cite{Isa}). They are never in rigid normal form and they are affinely equivalent to one of the following
\begin{align*}
v&=x^2\\
v&=e^x\\
\sin v&= e^x\\
e^v+e^x&=1
\end{align*}

This corresponds to the solutions for (\ref{reducedpde}) listed below
\begin{align*}
\tilde{f}&=0\\
\tilde{f}&=\frac12\\
\tilde{f}&=\tan x\\
\tilde{f}&=-\tanh x
\end{align*}

which integrates to 
\begin{align*}
f&=1\\
f&=x\\
f&=-2 \log \cos x\\
f&=-2 \log \cosh x.
\end{align*}

If follows
\begin{align*}
h&= \frac12 x^2 \\
h&= \e^x\\
h&= -\log \cos x \\
h&= \log \cosh x
\end{align*}
which yields the tubes 
\begin{align*}
v&=x^2\\
v&=\e^x\\
\e^v&= \cos x\\
\e^v&= \cosh x.
\end{align*}

\end{document}